\documentclass[11pt]{article}

\usepackage{amsthm,amsmath,amssymb,latexsym,lineno,cite}
\newtheorem{thm}{Theorem}[section]
\newtheorem{lem}[thm]{Lemma}
\newtheorem{cor}[thm]{Corollary}

\newtheorem{prop}[thm]{Proposition}

\theoremstyle{remark}

\renewcommand{\qed}{{\hfill\rule{4pt}{7pt}}\medskip}

\numberwithin{equation}{section}

\numberwithin{equation}{section}

\begin{document}
\begin{center}
{\Large\bf Multiple extensions of a finite Euler's pentagonal number theorem and the Lucas formulas}
\end{center}
\vskip 2mm \centerline{Victor J. W. Guo$^1$  and Jiang
Zeng$^{2}$}

\begin{center}
{\footnotesize $^1$Department of Mathematics, East China Normal
University, Shanghai 200062,
 People's Republic of China\\
{\tt jwguo@math.ecnu.edu.cn,\quad http://math.ecnu.edu.cn/\textasciitilde{jwguo}}\\[10pt]
$^2$ Universit\'e de Lyon; Universit\'e Lyon 1; Institut
Camille Jordan, UMR 5208 du CNRS ; 43, boulevard du 11 novembre 1918,
F-69622 Villeurbanne Cedex, France\\
{\tt zeng@math.univ-lyon1.fr,\quad
http://math.univ-lyon1.fr/\textasciitilde{zeng}} }
\end{center}

\vskip 0.7cm \noindent{\bf Abstract.} Motivated by the resemblance
of a multivariate series identity  and a finite analogue of
Euler's pentagonal number theorem, we  study  multiple extensions
of the latter formula. In a different direction we derive a common
extension of this multivariate series identity and two  formulas
of Lucas. Finally we give a combinatorial proof of Lucas' formulas.

\vskip 2mm
\noindent{\bf MR Subject Classifications}: 05A10; 11B65

\section{Introduction}
In a recent work \cite{GZ} we stumbled upon a multivariate identity involving
binomial coefficients (see \eqref{eq:dejavu}), which implies the
following identity:
\begin{align}
\sum_{r_1,\ldots,r_m\leq n}
\prod_{k=1}^{m}{n-r_k\choose r_{k+1}}
\frac{(-x)^{r_k}}
{(1+x)^{2r_k}}
=\frac{1-x^{m(n+1)}}
{(1-x^m)(1+x)^{mn}}, \label{eq:cyclic-x}
\end{align}
where $r_{m+1}=r_1.$ It is easy to see that the $x=\omega:=\frac{-1\pm i\sqrt{3}}{2}$
case of \eqref{eq:cyclic-x} reduces to
\begin{align}
\sum_{r_1,\ldots,r_m\leq n}
\prod_{k=1}^{m}{n-r_k\choose r_{k+1}}(-1)^{r_k}
=\begin{cases}
(-1)^{mn} (n+1),&\text{if $m\equiv 0 \pmod 3$},\\[5pt]
\displaystyle\frac{1-\omega^{m(n+1)}}
{(1-\omega^m)(1+\omega)^{mn}},&\text{if $m\not\equiv 0 \pmod 3$}.
\end{cases} \label{eq:cyc=-1}
\end{align}
This paper was motivated by the connection of \eqref{eq:cyc=-1} with some classical
formulas in the literature.

First of all, when  $m=1$, the formula
\eqref{eq:cyc=-1} has a known $q$-analogue (see \cite{BG,Cigler,EZ,Warnaar}) as follows:
\begin{align}
\sum_{k=0}^{\lfloor n/2\rfloor}(-1)^k q^{k\choose 2}{n-k\brack k}
=\begin{cases}
(-1)^{\lfloor n/3\rfloor} q^{n(n-1)/6},&\text{if $n\not\equiv 2 \pmod 3$},\\[5pt]
0,&\text{if $n\equiv 2 \pmod 3$}.
\end{cases} \label{eq:zeil}
\end{align}
where the {\it $q$-binomial coefficient ${n\brack k}$}
is defined by
$$
{n\brack k}
=\begin{cases}
 \displaystyle\prod_{i=1}^{k}\frac{1-q^{n-i+1}}{1-q^i}, & \text{if $0\leq k\leq n$,} \\
 0, &\text{otherwise.}
 \end{cases}
$$

Replacing $n$ by $3L$ or $3L+1$ and $q$ by $1/q$ in \eqref{eq:zeil} yields
\begin{align}
&\sum_{j=-L}^{L}(-1)^j q^{j(3j+1)/2}{2L-j\brack L+j}=1, \label{eq:finiteeuler-1}\\
&\sum_{j=-L}^{L}(-1)^j q^{j(3j-1)/2}{2L-j+1\brack L+j}=1, \label{eq:finiteeuler-2}
\end{align}
as mentioned in \cite{Warnaar}. Both \eqref{eq:finiteeuler-1} and \eqref{eq:finiteeuler-2} reduce to
Euler's pentagonal number theorem~\cite[p.~11]{And} when $L\to \infty$:
\begin{align} \label{eq:euler}
\sum_{j=-\infty}^{\infty}(-1)^jq^{j(3j-1)/2}=\prod_{n=1}^\infty(1-q^n).
\end{align}
It is then  natural to look for multiple analogues of
\eqref{eq:zeil} in light of \eqref{eq:cyc=-1}.
This will be the main object of Section 2.

Secondly, as will be shown, Eq.~\eqref{eq:cyclic-x} is also related
to the two formulas of Lucas (cf. \cite{Gould}):
\begin{align}
&\sum_{k=0}^{\lfloor n/2\rfloor}{n-k\choose k}(x+y)^{n-2k}(-xy)^k =\frac{x^{n+1}-y^{n+1}}{x-y},\label{lucas1}\\
&\sum_{k=0}^{\lfloor n/2\rfloor}\frac{n}{n-k}{n-k\choose k}(x+y)^{n-2k}(-xy)^k =x^n+y^n.\label{lucas2}
\end{align}
In section~3 by using the multivariate Lagrange inversion formula
we will prove a generalization of the formula~\eqref{eq:cyclic-x},
which is also a common  extension of Lucas' formulas
\eqref{lucas1} and \eqref{lucas2}.

Finally, as Shattuck and Wagner~\cite{SW} have recently given combinatorial a proof
of \eqref{lucas1} and \eqref{lucas2} with $x=1$ and $y=\omega$, we shall give a combinatorial
proof of Lucas' formulas in their full generality in Section~4.

We conclude this section with some remarks. It is known (see
\cite{Cigler}) that \eqref{eq:zeil} is actually equivalent to an
identity due to Rogers (see \cite[p.~29, Example 10]{And}). Some
modern proofs are given by Ekhad and Zeilberger~\cite{EZ} and
Warnaar~\cite{Warnaar}. The reader is also referred to Cigler's
paper~\cite{Cigler} for more information and proofs of
\eqref{eq:zeil}. Some known multiple and finite extensions of Euler's
pentagonal number theorem \eqref{eq:euler} can be found in
\cite{Andrews,Milas}, \cite[(6.2)]{GK}, \cite[(1)]{Krattenthaler}
and  the references therein. Note also that the $x+y=1$ and $xy=z$ cases of \eqref{lucas1} and
\eqref{lucas2} are sometimes called the Binet formulas~(see
\cite[p.~204]{GKP}).

\section{Common  extensions of   \eqref{eq:cyc=-1} and \eqref{eq:zeil} }
We shall adopt the standard notation of $q$-series in \cite{GR}. Let
$$
(a;q)_n=(1-a)(1-aq)\cdots(1-aq^{n-1}),\quad n=0,1,2,\ldots.
$$
Then the $q$-Chu-Vandermonde formula can be written as:
\begin{align}\label{eq:qchuvan}
 \sum_{k\geq 0}\frac{(a;q)_k(q^{-N};q)_k}{(c;q)_k(q;q)_k}\left(\frac{cq^N}{a}\right)^k
=\frac{(c/a;q)_N}{(c;q)_N}
\end{align}
(see \cite[p.~354]{GR}). We need the following two variations of \eqref{eq:qchuvan}.
\begin{lem} Let $n\geq 1$ and $r,t\leq n$. Then
\begin{align}
&\sum_{s=0}^{n-r}{n-r\brack s}{n-s\brack t}q^{s\choose 2}(-1)^s
=q^{(n-r)(n-t)}{r\brack n-t},  \label{eq:nrst-q} \\
&\sum_{s=0}^{n-r}{n-r\brack s}{n-s\brack t}q^{s(s+2r+2t-2n+1)/2}(-1)^s
={r\brack n-t}.  \label{eq:nrst-q2}
\end{align}
\end{lem}
Indeed, Eq.~\eqref{eq:nrst-q} follows from \eqref{eq:qchuvan}
with $a=q^{r-n}$, $N=n-t$ and $c=q^{-n}$, and
\eqref{eq:nrst-q2} can be derived from  \eqref{eq:nrst-q}
by the substitution $q\to q^{-1}$.

\begin{thm}Let $m,n\geq 1$ and $x_{3k}=-1$ for all $1\leq k\leq m$. Then
\begin{align}
\sum_{r_1,\ldots,r_{3m}\leq n}\prod_{k=1}^{3m}
{n-r_k\brack r_{k+1}}q^{r_k\choose 2}x_k^{r_k}
=\frac{(x_1x_4\cdots x_{3m-2})^{n+1}-(x_2x_5\cdots x_{3m-1})^{n+1}}
{x_1x_4\cdots x_{3m-2}-x_2x_5\cdots x_{3m-1}}q^{m{n\choose 2}},
\label{eq:multi-3m}
\end{align}
where $r_{3m+1}=r_1$.
\end{thm}

\begin{proof}
By \eqref{eq:nrst-q}, the left-hand side of \eqref{eq:multi-3m} equals
\begin{align}
&\hskip -3mm
\sum_{\substack{r_{3i-2},r_{3i-1}\leq n\\ 1\leq i\leq m} }
\prod_{k=1}^{m}
{n-r_{3k-2}\brack r_{3k-1}}q^{{r_{3k-2}\choose 2}+{r_{3k-1}\choose 2}}
x_{3k-2}^{r_{3k-2}}x_{3k-1}^{r_{3k-1}} \nonumber\\
&\quad\times
\sum_{r_3,r_6,\ldots,r_{3m}\leq n}
\prod_{k=1}^{m}{n-r_{3k-1}\brack r_{3k}}{n-r_{3k}\brack r_{3k+1}}
q^{r_{3k}\choose 2}(-1)^{r_{3k}}  \nonumber\\
&=\sum_{\substack{r_{3i-2},r_{3i-1}\leq n\\ 1\leq i\leq m} }
\prod_{k=1}^{m}
{n-r_{3k-2}\brack r_{3k-1}}{r_{3k-1}\brack n-r_{3k+1}}
q^{{r_{3k-2}\choose 2}+{r_{3k-1}\choose 2}+(n-r_{3k-1})(n-r_{3k+1})}
x_{3k-2}^{r_{3k-2}}x_{3k-1}^{r_{3k-1}} . \label{eq:multi-3k-12}
\end{align}
Note that
\begin{align*}
&\prod_{k=1}^{m}
{n-r_{3k-2}\brack r_{3k-1}}{r_{3k-1}\brack n-r_{3k+1}}\\
&=\begin{cases}
    1, & \hbox{if $r_{3k-2}+r_{3k-1}=n$ and $r_{3k-1}+r_{3k+1}=n$ for all $1\leq k\leq m$,} \\
    0, & \hbox{otherwise.}
\end{cases}
\end{align*}
 Therefore, the nonzero terms in the right-hand side of
\eqref{eq:multi-3k-12} are those indexed by $r_1=r_4=\cdots=r_{3m-2}$ and
$r_2=r_5=\cdots=r_{3m-1}=n-r_1$.  Finally, since
$$
{r_{3k-2}\choose 2}+{r_{3k-1}\choose 2}+(n-r_{3k-1})(n-r_{3k+1})={n\choose 2}
$$
for $r_{3k-2}+r_{3k-1}=n$ and $r_{3k-1}+r_{3k+1}=n$,  we see that the right-hand side of
\eqref{eq:multi-3k-12} equals
$$
\sum_{i=0}^n q^{m{n\choose 2}}(x_1x_4\cdots x_{3m-2})^{i}(x_2x_5\cdots x_{3m-1})^{n-i},
$$
as desired.
\end{proof}

Letting $x_{k}=-1$ for all $1\leq k\leq 3m$ in the above theorem
yields a $q$-analogue of \eqref{eq:cyc=-1} for $m\equiv 0\pmod{3}$.
\begin{cor} Let $m,n\geq 1$. Then
\begin{align}
\sum_{r_1,\ldots,r_{3m}\leq n}\prod_{k=1}^{3m}
{n-r_k\brack r_{k+1}}q^{r_k\choose 2}(-1)^{r_k}
=(-1)^{mn}(n+1)q^{m{n\choose 2}},
\label{eq:multi-3msp}
\end{align}
where $r_{3m+1}=r_1$.
\end{cor}
The following theorem gives a $q$-analogue of \eqref{eq:cyc=-1} for  $m\not\equiv 0\pmod 3$.
\begin{thm}\label{thm:multi-zeil}
Let $m,n\geq 1$ and $m\not\equiv 0\pmod 3$. Then
\begin{align}
\sum_{r_1,\ldots,r_{m}\leq n}\prod_{k=1}^{m}
{n-r_k\brack r_{k+1}}q^{r_k\choose 2}(-1)^{r_k}
=\begin{cases}
(-1)^{\lfloor (m+n-1)m/3\rfloor} q^{mn(n-1)/6},&\text{if $n\not\equiv 2 \pmod 3$},\\[5pt]
0,&\text{if $n\equiv 2 \pmod 3$},
\end{cases} \label{eq:multi-zeil}
\end{align}
where $r_{m+1}=r_1$.
\end{thm}
\begin{proof}
Replacing $q$ by $q^{-1}$ in \eqref{eq:zeil}, we get
\begin{align}
\sum_{k=0}^{\lfloor n/2\rfloor}(-1)^{n-k} q^{k^2+{n-k\choose 2}}{n-k\brack k}
=\begin{cases}
(-1)^{\lfloor (2n+2)/3\rfloor} q^{n(n-1)/3},&\text{if $n\not\equiv 2 \pmod 3$},\\[5pt]
0,&\text{if $n\equiv 2 \pmod 3$}.
\end{cases} \label{eq:zeilnew}
\end{align}
By \eqref{eq:nrst-q}, we have
\begin{align}
\sum_{r_1=0}^n \sum_{r_2=0}^n {n-r_1\brack r_2}{n-r_2\brack r_1}
q^{{r_1\choose 2}+{r_2\choose 2}}(-1)^{r_1+r_2}
&=\sum_{r_1=0}^n  {r_1\brack n-r_1}
q^{{r_1\choose 2}+(n-r_1)^2}(-1)^{r_1}, \label{eq:r1r2}
\end{align}
which is the left-hand side of \eqref{eq:zeilnew}.
This proves the $m=2$ case of \eqref{eq:multi-zeil}.

Again, by \eqref{eq:nrst-q},  we see that
\begin{align*}
\sum_{r_1,\ldots,r_{4}\leq n}\prod_{k=1}^{4}
{n-r_k\brack r_{k+1}}q^{r_k\choose 2}(-1)^{r_k}
&=\sum_{r_1=0}^n \sum_{r_3=0}^n {r_1\brack n-r_3}{r_3\brack n-r_1}
  q^{{r_1\choose 2}+{r_3\choose 2}+2(n-r_1)(n-r_3)}(-1)^{r_1+r_3} \\
&=\sum_{r_1=0}^n \sum_{r_3=0}^n {n-r_1\brack r_3}{n-r_3\brack r_1}
  q^{{n-r_1\choose 2}+{n-r_3\choose 2}+2r_1 r_3} (-1)^{r_1+r_3},
\end{align*}
where $r_5=r_1$, is the product of the $q\to q^{-1}$ case of the left-hand side of \eqref{eq:r1r2}
and $q^{n(n-1)}$. This proves the $m=4$ case of \eqref{eq:multi-zeil}.

For $m>4$, by \eqref{eq:nrst-q} and \eqref{eq:nrst-q2}, there holds
\begin{align}
\sum_{r_1,\ldots,r_{4}\leq n}\prod_{k=1}^{4}
{n-r_k\brack r_{k+1}}q^{r_k\choose 2}(-1)^{r_k}
&=\sum_{r_1=0}^n \sum_{r_3=0}^n {r_1\brack n-r_3}{r_3\brack n-r_5}
  q^{{r_1\choose 2}+{r_3\choose 2}+(2n-r_1-r_5)(n-r_3)}(-1)^{r_1+r_3} \nonumber \\
&=\sum_{r_1=0}^n \sum_{r_3=0}^n {n-r_1\brack r_3}{n-r_3\brack n-r_5}
  q^{{n-r_1\choose 2}+{n-r_3\choose 2}+(n+r_1-r_5)r_3} (-1)^{r_1+r_3} \nonumber \\
&=q^{n\choose 2}\sum_{r_1=0}^n {r_1\brack r_5} q^{{n-r_1\choose 2}} (-1)^{r_1} \nonumber \\
&=(-1)^n q^{n\choose 2} \sum_{r_1=0}^n {n-r_1\brack r_5} q^{{r_1\choose 2}} (-1)^{r_1}.
\label{eq:4to1}
\end{align}
It follows that
\begin{align*}
\sum_{r_1,\ldots,r_{m}\leq n}\prod_{k=1}^{m}
{n-r_k\brack r_{k+1}}q^{r_k\choose 2}(-1)^{r_k}
=(-1)^n q^{n\choose 2}\sum_{r_1,\ldots,r_{m-3}\leq n}\prod_{k=1}^{m-3}
{n-r_k\brack r_{k+1}}q^{r_k\choose 2}(-1)^{r_k}.
\end{align*}
By induction we can complete the proof based on the $m=2,4$ cases.
\end{proof}

The following result gives multiple extensions of \eqref{eq:finiteeuler-1} and \eqref{eq:finiteeuler-2}.
\begin{cor}Let $L,m\geq 1$. Then
\begin{align}
\sum_{j_1,\ldots,j_m=-L}^{2L}\prod_{k=1}^{m}
{2L-j_k\brack L+j_{k+1}}q^{j_k j_{k+1}+{j_k+1\choose 2}}(-1)^{j_k}
&=\begin{cases}
    1, &\hbox{if $m\not\equiv 0\pmod 3$,} \\
    3L+1, &\hbox{if $m\equiv 0\pmod 3$,}
\end{cases} \label{eq:3ell-0}\\[5pt]
\sum_{j_1,\ldots,j_m=-L}^{2L+1}\prod_{k=1}^{m}
{2L-j_k+1\brack L+j_{k+1}}q^{j_k j_{k+1}+{j_k\choose 2}}(-1)^{j_k}
&=\begin{cases}
    (-1)^{\lfloor m^2/3\rfloor}, & \hbox{if $m\not\equiv 0\pmod 3$,} \\
    (-1)^{m/3}(3L+2), & \hbox{if $m\equiv 0\pmod 3$.}
\end{cases}
\label{eq:3ell-1}
\end{align}
where $j_{m+1}=j_1$.
\end{cor}

\begin{proof} Take $n=3L$ in \eqref{eq:multi-3msp} and \eqref{eq:multi-zeil},
and replace $r_k$ by $j_k+L$ and $q$ by $1/q$.
After making some simplifications, we obtain \eqref{eq:3ell-0}. In much the same way,
when $n=3L+1$ we are led to \eqref{eq:3ell-1}.
\end{proof}

For $m\geq 4$, we can further generalize
Theorem \ref{thm:multi-zeil} as in the following two theorems.

\begin{thm}Let $m\geq 4$, $n\geq 1$ and $m\equiv 1\pmod 3$. Let $s\leq m$ be a
positive integer such that $s\not\equiv 0\pmod 3$.
Then
\begin{align}
\sum_{r_1,\ldots,r_{m}\leq n}
z^{r_1-r_s}\prod_{k=1}^{m}{n-r_k\brack r_{k+1}}q^{r_k\choose 2}(-1)^{r_k}
=\begin{cases}
(-1)^{\lfloor (m+n-1)m/3\rfloor} q^{mn(n-1)/6},&\text{if $n\not\equiv 2\pmod 3$},\\[5pt]
0,&\text{if $n\equiv 2 \pmod 3$},
\end{cases} \label{eq:zeil-ref}
\end{align}
where $r_{m+1}=r_1$.
\end{thm}
\begin{proof}
We first prove the $m=4$ case. By symmetry, we may assume that $s=2$. In this case, the left-hand
side of \eqref{eq:zeil-ref} equals
\begin{align}
&\hskip -3mm
\sum_{r_1,\ldots r_4\leq n}
z^{r_1-r_2}
{n-r_1\brack r_2}{n-r_2\brack r_3}{n-r_3\brack r_4}{n-r_4\brack r_1}
q^{{r_1\choose 2}+\cdots+{r_4\choose 2}} (-1)^{r_1+\cdots+r_4}  \nonumber \\
&=\sum_{k=-n}^n \sum_{r_2,r_3,r_4\leq n}
z^k {n-r_2-k\brack r_2}{n-r_2\brack r_3}{n-r_3\brack r_4}{n-r_4\brack r_2+k}
q^{{r_2+k\choose 2}+{r_2\choose 2}+{r_3\choose 2}+{r_4\choose 2}}(-1)^{k+r_3+r_4}.
\label{eq:zeil-m=4}
\end{align}

By \eqref{eq:nrst-q}, for $k>0$, we have
$$
\sum_{r_3\leq n}
{n-r_2\brack r_3}{n-r_3\brack r_4}{n-r_4\brack r_2+k}q^{r_3\choose 2}(-1)^{r_3}
={r_2\brack n-r_4}{n-r_4\brack r_2+k}q^{(n-r_2)(n-r_4)}=0,
$$
while for $k<0$, we have
$$
\sum_{r_4\leq n}
{n-r_2\brack r_3}{n-r_3\brack r_4}{n-r_4\brack r_2+k}q^{r_4\choose 2}(-1)^{r_4}
={n-r_2\brack r_3}{r_3\brack n-r_2-k}q^{(n-r_3)(n-r_2-k)}=0.
$$
Therefore, the right-hand side of \eqref{eq:zeil-m=4} is
independent of $z$. This completes the proof of
\eqref{eq:zeil-ref} for $m=4$.

For $m\geq 7$, again by symmetry, we may assume that $s\geq (m+3)/2\geq 5$.
We then complete the proof by induction on $m$ and using \eqref{eq:4to1}.
\end{proof}

\begin{thm}Let $m\geq 5$, $n\geq 1$ and $m\equiv 2\pmod 3$.
Let $s\leq m$ be a positive integer such that $s\not\equiv 2\pmod 3$.
Then
\begin{align}
\sum_{r_1,\ldots,r_{m}\leq n}
z^{r_1-r_s}\prod_{k=1}^{m}{n-r_k\brack r_{k+1}}q^{r_k\choose 2}(-1)^{r_k}
=\begin{cases}
(-1)^{\lfloor (m+n-1)m/3\rfloor} q^{mn(n-1)/6},&\text{if $n\not\equiv 2\pmod 3$},\\[5pt]
0,&\text{if $n\equiv 2 \pmod 3$},
\end{cases}  \label{eq:zeil-ref2}
\end{align}
where $r_{m+1}=r_1$.
\end{thm}
\begin{proof}
For $m=5$, by symmetry, we may assume that $s=3$. In this case, the left-hand side of \eqref{eq:zeil-ref2}
may be written as
\begin{align}
&\sum_{k=-n}^n \sum_{r_2,\ldots,r_5\leq n}
z^k {n-r_3-k\brack r_2}{n-r_2\brack r_3}{n-r_3\brack r_4}{n-r_4\brack r_5}{n-r_5\brack r_3+k} \nonumber\\
&\quad\times
q^{{r_3+k\choose 2}+{r_2\choose 2}+{r_3\choose 2}+{r_4\choose 2}+{r_5\choose 2}}(-1)^{k+r_2+r_4+r_5}.
\label{eq:zeil-m=5}
\end{align}

By \eqref{eq:nrst-q}, for $k>0$, we have
$$
\sum_{r_4\leq n}
{n-r_3\brack r_4}{n-r_4\brack r_5}{n-r_5\brack r_3+k}q^{r_4\choose 2}(-1)^{r_4}
=0,
$$
while for $k<0$, we have
$$
\sum_{r_5\leq n}
{n-r_3\brack r_4}{n-r_4\brack r_5}{n-r_5\brack r_3+k}q^{r_5\choose 2}(-1)^{r_5}
=0.
$$
Therefore, the right-hand side of \eqref{eq:zeil-m=5} is independent of $z$. This completes
the proof of the $m=5$ case of \eqref{eq:zeil-ref2}.

For $m\geq 8$, again by symmetry, we may assume that $s\geq (m+3)/2$.
We then complete the proof by induction on $m$ and using \eqref{eq:4to1}.
\end{proof}

\section{Generalization of \eqref{eq:cyclic-x} and Lucas' formulas}
The following identity \eqref{eq:dejavu} was  already announced in
\cite{GZ}.
\begin{thm}
We have
\begin{align}\label{eq:dejavu}
\sum_{r_1,\ldots,r_m\leq n}
\prod_{k=1}^{m}{n-r_k\choose r_{k+1}}
\frac{(-x_k)^{r_k}}
{(1+x_k)^{r_k+r_{k+1}}}
=\frac{1-x_1^{n+1}\cdots x_m^{n+1}}
{1-x_1\cdots x_m}\prod_{k=1}^{m}\frac{1}{(1+x_k)^n},
\end{align}
where $r_{m+1}=r_1$.
\end{thm}
To prove this theorem, we need the following form of the multivariate
Lagrange inversion formula (see \cite[p.~21]{GJ}).
\begin{lem}\label{lem:lagrange}
Let $m\geq 1$ be a positive integer and ${\mathbf x}=(x_1,\ldots,x_m)$.
Suppose that  $x_i=u_i\phi_i({\bf x})$  for $i=1,\ldots,m$ and
$\phi_i$ is a formal power series in ${\bf x}$
 with complex coefficients such that
$\phi_i(0, \ldots, 0)\neq 0$. Then any formal power series $f({\mathbf x})$
with complex coefficients can be expanded into a power
series in ${\mathbf u}=(u_1,\ldots,u_m)$ as follows:
\[
f({\bf x}({\bf u}))=\sum_{\mathbf r\in \mathbb{N}^m}{\bf u}^{\bf r}
[{\bf x}^{\bf r}]
\left\{f({\bf x})\phi_1^{r_1}({\bf x})\ldots \phi_m^{r_m}({\bf x})
\Delta_m\right\},
\]
where $[{\bf x}^{\bf r}]f({\bf x})$ denotes the coefficient of
${\mathbf  x}^{\mathbf r}=x_1^{r_1}\ldots x_m^{r_m}$ in the series $f({\bf x})$ and
$$
\Delta_m=\det\left(\delta_{ij}-\frac{x_j}{\phi_i({\bf x})}
\frac{\partial \phi_i({\bf x})}{\partial x_j}\right)_{1\leq i,j\leq m}.
$$
\end{lem}
\noindent{\it Proof of \eqref{eq:dejavu}.}
Let $\phi_i({\mathbf x})=(1+x_{i-1})(1+x_i)$ ($1\leq i\leq m$),
where $x_{0}=x_m$.
Then $\Delta_1=\frac{1-x_1}{1+x_1}$ and for $m\geq 2$
\begin{align*}
\Delta_m &=\left|
\begin{array}{ccccc}
\displaystyle \frac{1}{1+x_1}    & 0  &\cdots & 0 &\displaystyle \frac{-x_m}{1+x_m}  \\[12pt]
\displaystyle \frac{-x_1}{1+x_1}  & \displaystyle \frac{1}{1+x_2}  & 0  &\cdots & 0   \\[12pt]
0 & \displaystyle \frac{-x_2}{1+x_2}& \displaystyle \frac{1}{1+x_3} &\cdots &    0    \\[12pt]
\vdots & \vdots & \vdots  & \vdots & \vdots  \\[12pt]
0 & \cdots &0 & \displaystyle \frac{-x_{m-1}}{1+x_{m-1}}& \displaystyle \frac{1}{1+x_m}\\[12pt]
\end{array}\right|
=\frac{1-x_1\cdots x_m}{\prod_{k=1}^{m}(1+x_k)}.
\end{align*}
Now take
\[
f({\mathbf x})
=\frac{1-x_1^{n+1}\cdots x_m^{n+1}}
{1-x_1\cdots x_m}\prod_{k=1}^{m}\frac{1}{(1+x_k)^{n}}.
\]
Then
\begin{align*}
f({\bf x})\phi_1^{r_1}({\bf x})\ldots \phi_m^{r_m}({\bf x})\Delta_m
=\frac{1-x_1^{n+1}\cdots x_m^{n+1}}
{\prod_{k=1}^{m}(1+x_k)^{n+1-r_k-r_{k+1}}}.
\end{align*}

Note that
\begin{align*}
[{\mathbf x}^{\mathbf r}]
\prod_{k=1}^{m}\frac{1}{(1+x_k)^{n+1-r_k-r_{k+1}}}
=\prod_{k=1}^{m}(-1)^{r_k}{n-r_{k+1}\choose r_k}
=\prod_{k=1}^{m}(-1)^{r_k}{n-r_{k}\choose r_{k+1}}.
\end{align*}
Also,
\begin{align*}
[{\mathbf x}^{\mathbf r}]
\prod_{k=1}^{m}\frac{x_k^{n+1}}{(1+x_k)^{n+1-r_k-r_{k+1}}}
=\begin{cases}
    \displaystyle\prod_{k=1}^{m}(-1)^{r_k}{n-r_k\choose r_{k+1}},
    &\hbox{if $r_1,\ldots,r_m\geq n+1$,} \\[5pt]
    0, & \hbox{otherwise.}
\end{cases}
\end{align*}
By subtraction we derive  from Lemma 4.1 that
\begin{align*}
f({\mathbf x})
&=\sum_{\min\{r_1,\ldots,r_m\}\leq n}
u_1^{r_1}\cdots u_m^{r_m}\prod_{k=1}^{m}(-1)^{r_k}{n-r_k\choose r_{k+1}} \\[5pt]
&=\sum_{r_1,\ldots,r_m\leq n}
\prod_{k=1}^{m}{n-r_k\choose r_{k+1}}
\frac{(-x_k)^{r_k}}{(1+x_k)^{r_k+r_{k+1}}},
\end{align*}
as desired. \qed

\noindent{\it Remark}. Strehl~\cite{St} has obtained more binomial
coefficients formulas by applying the multivariate Lagrange inversion formula.

Letting $x_i=x$ for all $i$ in \eqref{eq:dejavu} we obtain
\eqref{eq:cyclic-x}, while  letting $x=\frac{\sqrt{5}-3}{2}$ in \eqref{eq:cyclic-x}
we obtain the following remarkable identity
\begin{prop}\label{cor:cyclic-x}
For $m,n\geq 1$, we have
\begin{align*}
&\sum_{r_1,\ldots,r_m\leq n}
\prod_{k=1}^{m}{n-r_k\choose r_{k+1}}
=\frac{2^{m(n+1)}-\left(\sqrt{5}-3\right)^{m(n+1)}}
{\left(2^m-\left(\sqrt{5}-3\right)^m\right)\left(\sqrt{5}-1\right)^{mn}}, 
\end{align*}
where $r_{m+1}=r_1$.
\end{prop}

To see that \eqref{eq:cyclic-x} is a common multiple extension of
two formulas of Lucas, we first recall the following elementary
counting results (see, for example, \cite[Lemma 2.3.4]{Stanley}).
\begin{lem}\label{cor:key}
The number of ways of choosing $k$ points,
no two consecutive, from a collection of $n-1$ points arranged on a line is
${n-k\choose k}$. The number of ways of choosing $k$ points,
no two consecutive, from a collection of $n$ points arranged on a cycle is
$\frac{n}{n-k}{n-k\choose k}$.
\end{lem}

Now, the $m=1$ case of  \eqref{eq:cyclic-x} corresponds to
\begin{align}\label{eq:sqrt-fibo}
\sum_{k=0}^{\lfloor n/2\rfloor}{n-k\choose k}\frac{(-x)^k}{(1+x)^{2k}}
=\frac{1-x^{n+1}}{(1-x)(1+x)^n}.
\end{align}
On the other hand, for $r_1,\ldots, r_{m}\in \{0,1\}$ and $r_{m+1}=r_1$,
the product $\prod_{k=1}^{m}{1-r_k\choose r_{k+1}}$  equals $1$ if there are no two
consecutive $1$'s in the sequence $r_1,\ldots,r_m,r_{m+1}$, and
$0$ otherwise. Thus, by Lemma~\ref{cor:key},
 the $n=1$ case of
\eqref{eq:cyclic-x} corresponds
 to the following identity:
\begin{align}\label{eq:non-adjx}
\sum_{k=0}^{\lfloor m/2\rfloor}\frac{m}{m-k}{m-k\choose
k}\frac{(-x)^k}{(1+x)^{2k}} =\frac{1+x^m}{(1+x)^m}.
\end{align}
Clearly Lucas' formulas \eqref{lucas1} and \eqref{lucas2} are
equivalent to \eqref{eq:sqrt-fibo} and \eqref{eq:non-adjx}. When
$x=\omega$ the latter formulas  (replacing $m$ by $n$ in
\eqref{eq:non-adjx}) can be written as
\begin{align}
& \sum_{k=0}^{\lfloor n/2\rfloor}{n-k\choose k}(-1)^k=
\frac{1-\omega^{n+1}}{(1-\omega)(1+\omega)^n}
=\begin{cases}
1,&\text{if $n\equiv 0,1 \pmod 6$},\\
0,&\text{if $n\equiv 2,5 \pmod 6$},\\
-1,&\text{if $n\equiv 3,4 \pmod 6$},
\end{cases} \label{eq:6case2}
\end{align}
and
\begin{align}
\sum_{k=0}^{\lfloor n/2\rfloor}\frac{n}{n-k}{n-k\choose k}(-1)^k
=\frac{1+\omega^n}{(1+\omega)^n}=\begin{cases}
2,&\text{if $n\equiv 0 \pmod 6$},\\
1,&\text{if $n\equiv 1,5 \pmod 6$},\\
-1,&\text{if $n\equiv 2,4 \pmod 6$},\\
-2,&\text{if $n\equiv 3 \pmod 6$}.
\end{cases} \label{eq:6case}
\end{align}
Motivated by the recent combinatorial proof of \eqref{eq:6case2} and \eqref{eq:6case}
by Shattuck and Wagner \cite{SW}, we shall give a combinatorial proof of
a polynomial version of \eqref{eq:sqrt-fibo} and \eqref{eq:non-adjx} in the next section.

\section{Combinatorial proof of Lucas' formulas}
Letting  $m=\frac{-x}{1+x}$ in \eqref{eq:sqrt-fibo} and
\eqref{eq:non-adjx},
 we obtain
\begin{align}
&\sum_{k=0}^{\lfloor n/2\rfloor}{n-k\choose k}m^k(m+1)^k
=\frac{1}{2m+1}\left((m+1)^{n+1}-(-m)^{n+1}\right),\label{eq:bij-1} \\
&\sum_{k=0}^{\lfloor n/2\rfloor}\frac{n}{n-k}{n-k\choose
k}m^k(m+1)^k =(m+1)^{n}+(-m)^{n}. \label{eq:bij-2}
\end{align}
We now give a bijective proof of \eqref{eq:bij-1} and \eqref{eq:bij-2} assuming that $m$
is a positive integer. Obviously this is sufficient to prove their validity.

\begin{itemize}
\item For any positive integer $n$, let $[n]:=\{1,\ldots, n\}$.
Given $n>1$, let $\mathcal S$ be the set of all triples $(A;f,g)$ such that
$A$ is a subset of $[n-1]$ without consecutive integers, $f\colon
A\to [m]$ and $g\colon A\to [m+1]$ are two mappings (or
colorings). By Lemma~\ref{cor:key} the left-hand side of
\eqref{eq:bij-1} is the cardinality of $\mathcal S$.

A \emph{chain} is a set of consecutive integers, the cardinality being called its length.
Let $X$ be a set of integers. A chain $Y\subseteq X$ is called
\emph{maximal} if there is no other chain $Y'$ in $X$ such that $Y\subset Y'$.
It is clear that $X$ can be decomposed uniquely as a union of its
disjoint maximal chains.
Let $\mathcal T$ be the set of all pairs $(X;h)$
where $X\subseteq [n]$ such that the maximal chain containing $n$
in $X$ (if exists) is of \emph{even} length and $h\colon X\to [m]$ is a mapping.
Since the number of all pairs $(X;h)$ with $X\subseteq [n]$
and $h\colon X\to [m]$ is equal to $(m+1)^n$,  the  number of all such
pairs $(X;h)$ with the maximal chain containing $n$ being of even length, say $2k$, is given by
$$
m^{2k}(m+1)^{n-2k-1}=m^{2k}(m+1)^{n-2k}-m^{2k+1}(m+1)^{n-2k-1}
$$
if $2k<n$, and $m^n$ if $2k= n$. Summing up, the cardinality of $\mathcal T$ equals
$$
\sum_{k=0}^{\lfloor n/2\rfloor}m^{2k}(m+1)^{n-2k}-
\sum_{k=0}^{\lfloor
(n-1)/2\rfloor}m^{2k+1}(m+1)^{n-2k-1}=\sum_{k=0}^n(-m)^k(m+1)^{n-k},
$$
i.e., the right-hand side of \eqref{eq:bij-1}.

It remains to establish a bijection $\theta\colon \mathcal S\to\mathcal T$. For each
$(A;f,g)\in\mathcal S$, let $B=\{i+1\colon i\in A\text{ and }g(i)\in [m]\}$ and define
 $\theta (A;f,g)=(X;h)$ by $X=A\cup B$
and  $h|_A=f$ and $h(i)=g(i-1)$ for $i\in B$. It is easy to see
that $(X;h)\in\mathcal T$. Conversely, let $(X;h)\in\mathcal T$, suppose
$X=X_1\cup\cdots \cup X_s$, where $X_i$ is a maximal chain of $X$
for each $i=1,\ldots,s$. Write $X_i=\{x_{i,1},x_{i,2},x_{i,3},\ldots\}$ in increasing order.
Define the tripe $(A;f,g)\in\mathcal S$ by $A=\cup_{i=1}^s
\{x_{i,1},x_{i,3},x_{i,5},\ldots\}$, $f=h|_A$ and $g(i)=h(i+1)$ if
$i+1\in X\setminus A$ and $g(i)=m+1$ if $i+1\notin X\setminus A$.
Then $(A;f,g)$ is the unique preimage of $(X;h)$ under the mapping $\theta$.
This completes the proof of \eqref{eq:bij-1}.

\item Next consider the cyclic group $\mathbb{Z}_n=\{0,1,\ldots,n-1\}$.
Let $\mathcal U$ be the set of triples $(A;f,g)$, where $A$ is a subset of
$\mathbb{Z}_n$ without consecutive elements of $\mathbb{Z}_n$,
$f\colon A\to [m]$ and $g\colon A\to [m+1]$ are two mappings. By
Lemma~\ref{cor:key} the left-hand side of \eqref{eq:bij-2} is
equal to the cardinality of $\mathcal U$.

Let $\mathcal V$ be the set of all pairs $(X;h)$ where $X\subseteq \mathbb{Z}_n$ and
$h\colon X\to [m]$ is a mapping.
We define a mapping $\varphi\colon\mathcal U\to\mathcal V$
as follows.

For each $(A;f,g)\in\mathcal U$,  let $B=\{i+1\colon i\in A\text{ and }g(i)\in [m]\}$, $X=A\cup B$, $h|_A=f$ and
$h(i)=g(i-1)$ for $i\in B$. Then $\varphi(A;f,g)=(X;h)\in\mathcal V$. Conversely,
each $(X;h)\in\mathcal V$ with $X\subsetneq \mathbb{Z}_n$ has a unique preimage under the mapping $\varphi$.
However each $(\mathbb{Z}_n;h)\in\mathcal V$ has no preimage if $n$ is odd, and has
two preimages if $n$ is even: $(A_1;f_1,g_1)$ and $(A_2;f_2,g_2)$,
where $A_1=\{0,2,4,\ldots,n-2\}$, $A_2=\{1,3,5,\ldots,n-1\}$,
$f_1(i)=h(i)$ and $g_1(i)=h(i+1)$ for $i\in A_1$; $f_2(i)=h(i)$
and $g_2(i)=h(i+1)$ for $i\in A_2$. Thus, the cardinality of $\mathcal U$
is equal to $(m+1)^n+(-m)^{n}$. This completes the proof of
\eqref{eq:bij-2}.
\end{itemize}

\section*{Acknowledgments} This work was done during the first
author's visit to Institut Camille Jordan, Universit\'e
Claude Bernard (Lyon I), and was supported by a French postdoctoral fellowship.
The authors thank the two anonymous referees for
helpful comments on a previous version of this paper.

\renewcommand{\baselinestretch}{1}

\end{document}